\documentclass[11pt,a4paper]{amsart}

\usepackage{amssymb}
\usepackage[arrow,curve,matrix]{xy}

\title[Derived equivalences of triangular matrix rings]
{Derived equivalences of triangular matrix rings arising from
extensions of tilting modules}
\author{Sefi Ladkani}
\address{Einstein Institute of Mathematics, The Hebrew University of Jerusalem, Jerusalem 91904, Israel}
\email{sefil@math.huji.ac.il}

\newcommand{\iu}{i^{-1}}
\newcommand{\ius}{i^{!}}
\newcommand{\il}{i_{*}}
\newcommand{\ils}{i_{!}}
\newcommand{\ju}{j^{-1}}
\newcommand{\jus}{j^{\natural}}
\newcommand{\jl}{j_{!}}
\newcommand{\jls}{j_{*}}
\newcommand{\ii}{\il\iu}
\newcommand{\jj}{\jl\ju}

\newcommand{\cA}{\mathcal{A}}
\newcommand{\cB}{\mathcal{B}}
\newcommand{\cC}{\mathcal{C}}
\newcommand{\cD}{\mathcal{D}}

\newcommand{\gL}{\Lambda}

\newcommand{\dA}{\cD^b(\cA)}
\newcommand{\dB}{\cD^b(\cB)}
\newcommand{\dC}{\cD^b(\cC)}
\newcommand{\dL}{\cD^b(\gL)}

\newcommand{\bZ}{\mathbb{Z}}

\DeclareMathOperator{\Mod}{Mod} \DeclareMathOperator{\modf}{mod}
\DeclareMathOperator{\per}{per}

\DeclareMathOperator{\coker}{coker}

\DeclareMathOperator{\gldim}{gl.\!dim}
\DeclareMathOperator{\Hom}{Hom}
\DeclareMathOperator{\Ext}{Ext}
\DeclareMathOperator{\End}{End}
\DeclareMathOperator{\Id}{Id}
\DeclareMathOperator{\ob}{ob}
\DeclareMathOperator{\add}{add}

\DeclareMathOperator{\GL}{GL}

\theoremstyle{plain}
\newtheorem{thm}{Theorem}[section]
\newtheorem*{thm*}{Theorem}
\newtheorem{lem}[thm]{Lemma}
\newtheorem{prop}[thm]{Proposition}
\newtheorem{cor}[thm]{Corollary}
\newtheorem*{cor*}{Corollary}

\theoremstyle{definition}
\newtheorem{defn}[thm]{Definition}
\newtheorem{rem}[thm]{Remark}

\newtheorem{exmp}[thm]{Example}

\numberwithin{equation}{section}

\begin{document}

\begin{abstract}
A triangular matrix ring $\gL$ is defined by a triplet $(R,S,M)$ where
$R$ and $S$ are rings and $_RM_S$ is an $S$-$R$-bimodule. In the main
theorem of this paper we show that if $T_S$ is a tilting $S$-module,
then under certain homological conditions on the $S$-module $M_S$, one
can extend $T_S$ to a tilting complex over $\gL$ inducing a derived
equivalence between $\gL$ and another triangular matrix ring specified
by $(S', R, M')$, where the ring $S'$ and the $R$-$S'$-bimodule $M'$
depend only on $M$ and $T_S$, and $S'$ is derived equivalent to $S$.
Note that no conditions on the ring $R$ are needed.

These conditions are satisfied when $S$ is an Artin algebra of finite
global dimension and $M_S$ is finitely generated. In this case,
$(S',R,M') = (S, R, DM)$ where $D$ is the duality on the category of
finitely generated $S$-modules. They are also satisfied when $S$ is
arbitrary, $M_S$ has a finite projective resolution and $\Ext^n_S(M_S,
S) = 0$ for all $n > 0$. In this case, $(S',R,M') = (S, R, \Hom_S(M,
S))$.
\end{abstract}

\maketitle

\section{Introduction}
Triangular matrix rings and their homological properties have been
widely studied, see for example
\cite{Chase61,Fields70,FGR75,MichelenaPlatzeck00,PalmerRoos73}. The
question of derived equivalences between different such rings was
explored in the special case of one-point extensions of
algebras~\cite{BarotLenzing03}. Another aspect of this question was
addressed by considering examples of triangular matrix algebras of a
simple form, such as incidence algebras of posets~\cite{Ladkani08}. In
this paper we extend the results of~\cite{Ladkani08} to general
triangular matrix rings.

A triangular matrix ring $\gL$ is defined by a triplet $(R,S,M)$ where
$R$ and $S$ are rings and $_RM_S$ is an $S$-$R$-bimodule. The category
of (right) $\gL$-modules can be viewed as a certain gluing of the
categories of $R$-modules and $S$-modules, specified by four exact
functors. This gluing naturally extends to the bounded derived
categories. We note the similarity to the classical ``recollement''
situation, introduced by Beilinson, Bernstein and Deligne~\cite{BBD82},
involving six functors between three triangulated categories,
originally inspired by considering derived categories of sheaves on
topological spaces, and later studied for derived categories of modules
by Cline, Parshall and Scott~\cite{CPS88a,CPS88b}, see
also~\cite{Koenig91}.

In Section~\ref{sec:cat} we show that triangular matrix rings arise
naturally as endomorphism rings of certain rigid complexes over abelian
categories that are glued from two simpler ones. Here, a complex $T \in
\dC$ is \emph{rigid} if $\Hom_{\dC}(T,T[n])=0$ for all $n \neq 0$,
where $\dC$ denotes the bounded derived category of an abelian category
$\cC$. Similarly, an object $T \in \cC$ is \emph{rigid} if it is rigid
as a complex.

Indeed, when $\cC$ is glued from the abelian categories $\cA$ and
$\cB$, we construct, for any projective object of $\cA$ and a rigid
object of $\cB$ satisfying some homological conditions, a new rigid
complex in $\dC$ whose endomorphism ring is a triangular matrix ring.

In particular, as demonstrated in Section~\ref{sec:comma}, this
construction applies for comma categories defined by two abelian
categories $\cA$, $\cB$ and a right exact functor $F : \cA \to \cB$. In
this case, any projective $P$ of $\cA$ and a rigid object $T_{\cB} \in
\cB$ satisfying $\Ext^n_{\cB}(FP, T_{\cB}) = 0$ for all $n
> 0$, give rise to a rigid complex $T$ over the comma category, whose
endomorphism ring is a triangular matrix ring which can be explicitly
computed in terms of $P$, $T_{\cB}$ and $FP$, as
\[
\End_{\dC}(T) \simeq
\begin{pmatrix}
\End_{\cB}(T_{\cB}) & \Hom_{\cB}(FP, T_{\cB}) \\
0 & \End_{\cA}(P)
\end{pmatrix}.
\]

In Section~\ref{sec:alg} we apply this construction for categories of
modules over triangular matrix rings. For a ring $\gL$, denote by $\Mod
\gL$ the category of all right $\gL$-modules, and by $\dL$ its bounded
derived category. Recall that a complex $T \in \dL$ is a \emph{tilting
complex} if it is rigid and moreover, the smallest full triangulated
subcategory of $\dL$ containing $T$ and closed under forming direct
summands, equals $\per \gL$, the full subcategory in $\dL$ of complexes
quasi-isomorphic to \emph{perfect complexes}, that is, bounded
complexes of finitely generated projective $\gL$-modules. If, in
addition, $H^n(T) = 0$ for all $n \neq 0$, we call $T$ a \emph{tilting
module} and identify it with the module $H^0(T)$.

Two rings $\gL$ and $\gL'$ are \emph{derived equivalent} if $\dL$ and
$\cD^b(\gL')$ are equivalent as triangulated categories. By Rickard's
Morita theory for derived equivalence~\cite{Rickard89}, this is
equivalent to the existence of a tilting complex $T \in \dL$ such that
$\End_{\dL}(T) = \gL'$.

When $\gL$ is a triangular matrix ring defined by two rings $R$, $S$
and a bimodule $_RM_S$, the category $\Mod \gL$ is the comma category
of $\Mod R$, $\Mod S$ with respect to the functor $- \otimes M : \Mod R
\to \Mod S$. In this case, starting with the projective $R$-module $R$
and a tilting $S$-module $T_S$, the complex $T$ constructed in
Section~\ref{sec:comma} is not only rigid, but also a tilting complex,
hence we deduce a derived equivalence between $\gL$ and the triangular
matrix ring $\End_{\dL}(T)$, as expressed in the theorem below.

\begin{thm*}
Let $R, S$ be rings and $T_S$ a tilting $S$-module. Let $_RM_S$ be an
$S$-$R$-bimodule such that as an $S$-module, $M_S \in \per S$ and
$\Ext^n_S(M_S, T_S) = 0$ for all $n > 0$. Then the triangular matrix
rings
\[
\gL =
\begin{pmatrix}
R & M \\
0 & S \\
\end{pmatrix}
\quad \text{and} \quad \widetilde{\gL} =
\begin{pmatrix}
\End_S(T_S) & \Hom_S(M, T_S) \\
0 & R \\
\end{pmatrix}
\]
are derived equivalent.
\end{thm*}

We note that the assumption that $T_S$ is a tilting module implies that
the rings $S$ and $\End_S(T_S)$ are derived equivalent, hence the
triangular matrix ring specified by the triplet $(R, S, M)$ is derived
equivalent to a one specified by $(S', R, M')$ where $S'$ is derived
equivalent to $S$. We note also that no conditions on the ring $R$ (or
on $M$ as a left $R$-module) are necessary.

The above theorem has two interesting corollaries, corresponding to the
cases where $T_S$ is injective or projective.

For the first corollary, let $S$ be an Artin algebra, and let $D: \modf
S \to \modf S^{op}$ denote the duality. When $S$ has finite global
dimension, one can take $T_S$ to be the module $DS$ which is then an
injective tilting module.

\begin{cor*}
Let $R$ be a ring, $S$ an Artin algebra with $\gldim S < \infty$ and
$_RM_S$ an $S$-$R$-bimodule which is finitely generated as an
$S$-module. Then the triangular matrix rings
\[
\gL =
\begin{pmatrix}
R & M \\
0 & S \\
\end{pmatrix}
\quad \text{and} \quad \widetilde{\gL} =
\begin{pmatrix}
S & DM \\
0 & R \\
\end{pmatrix}
\]
are derived equivalent, where $D$ is the duality on $\modf S$.
\end{cor*}

The ring $\widetilde{\gL}$ depends only on $R$, $S$ and $M$, hence it
may be considered as a derived equivalent \emph{mate} of $\gL$.

\medskip

The second corollary of the above theorem is obtained by taking the
tilting $S$-module to be $S$.

\begin{cor*}
Let $R$, $S$ be rings and $_RM_S$ an $S$-$R$-bimodule such that as an
$S$-module, $M_S \in \per S$ and $\Ext^n_S(M_S, S) = 0$ for all $n >
0$. Then the triangular matrix rings
\[
\gL =
\begin{pmatrix}
R & M \\
0 & S \\
\end{pmatrix}
\quad \text{and} \quad \widetilde{\gL} =
\begin{pmatrix}
S & \Hom_S(M, S) \\
0 & R \\
\end{pmatrix}
\]
are derived equivalent.
\end{cor*}

This corollary applies to the following situations, listed in
descending order of generality; the ring $S$ is \emph{self-injective}
(that is, $S$ is injective as a module over itself) and $M_S$ is
finitely generated projective; the ring $S$ is semi-simple and $M_S$ is
finitely generated; the ring $S$ is a division ring and $M$ is finite
dimensional over $S$. The latter case implies that a triangular matrix
ring which is a one-point extension is derived equivalent to the
corresponding one-point co-extension.

In Section~\ref{sec:remarks} we conclude with three remarks concerning
the specific case of finite dimensional triangular matrix algebras over
a field.


First, in the case when $R$ and $S$ are finite dimensional algebras
over a field and both have finite global dimension, an alternative
approach to show the derived equivalence of the triangular matrix
algebras specified by $(R,S,M)$ and its mate $(S,R,DM)$ is to prove
that the corresponding repetitive algebras are isomorphic and then use
Happel's Theorem~\cite[(II,4.9)]{Happel88}. However, in the case that
only one of $R$ and $S$ has finite global dimension, Happel's Theorem
cannot be used, but the derived equivalence still holds. Moreover, as
we show in Example~\ref{ex:matenoneq}, there are cases when none of $R$
and $S$ have finite global dimension and the corresponding algebras are
not derived equivalent, despite the isomorphism between their
repetitive algebras.

Second, one can directly prove, using only matrix calculations, that
when at least one of $R$ and $S$ has finite global dimension, the
Cartan matrices of the triangular matrix algebra $(R,S,M)$ and its mate
are equivalent over $\bZ$, a result which is a direct consequence of
Theorem~\ref{t:RSMartin}.

Third, we note that in contrast to triangular matrix algebras, in the
more general situation of trivial extensions of algebras, the mates $A
\ltimes M$ and $A \ltimes DM$ for an algebra $A$ and a bimodule
$_AM_A$, are typically not derived equivalent.

\subsection*{Acknowledgement}
I am grateful to B.~Keller for discussions of a preliminary version
of~\cite{Ladkani08} which led to this current paper.

\section{The gluing construction}
\label{sec:cat} \setcounter{subsection}{1}

Let $\cA$, $\cB$, $\cC$ be three abelian categories. Similarly
to~\cite[(1.4)]{BBD82}, we view $\cC$ as glued from $\cA$ and $\cB$ if
there exist certain functors $\iu, \il, \ju, \jl$ as described below.
Note, however, that we start by working at the level of the abelian
categories and not their derived categories. In addition, the
requirement in~\cite{BBD82} of the existence of the additional adjoint
functors $\ius, \jls$ is replaced by the orthogonality
condition~\eqref{e:ijC0}.

\begin{defn}
A quadruple of additive functors $\iu, \il, \ju, \jl$ as in the diagram
\[
\xymatrix{\cA \ar@<1ex>[r]^{\il} & \cC \ar@<1ex>[l]^{\iu}
\ar@<1ex>[r]^{\ju} & \cB \ar@<1ex>[l]^{\jl}}
\]
is called \emph{gluing data} if it satisfies the four
properties~(\ref{ss:adj})--(\ref{ss:ortho}) below.
\end{defn}

\subsubsection{Adjunction}
\label{ss:adj}

$\iu$ is a left adjoint of $\il$ and $\ju$ is a right adjoint of $\jl$.
That is, there exist bi-functorial isomorphisms
\begin{align}
\label{e:iadj}
\Hom_{\cA}(\iu C, A) &\simeq \Hom_{\cC}(C, \il A) \\
\label{e:jadj} \Hom_{\cB}(B, \ju C) &\simeq \Hom_{\cC}(\jl B, C)
\end{align}
for all $A \in \ob \cA$, $B \in \ob \cB$, $C \in \ob \cC$.

\subsubsection{Exactness}
\label{ss:exact}

The functors $\iu, \il, \ju, \jl$ are exact.

Note that by the adjunctions above, we automatically get that the
functors $\il, \ju$ are left exact while $\iu, \jl$ are right exact.

\subsubsection{Extension}
\label{ss:ext}

For every $C \in \ob \cC$, the adjunction morphisms $\jj C \to C$ and
$C \to \ii C$ give rise to a short exact sequence
\begin{equation}
\label{e:jiSES} 0 \to \jj C \to C \to \ii C \to 0
\end{equation}

\subsubsection{Orthogonality}
\label{ss:ortho}

\begin{align}
\label{e:ij0} & \iu \jl = 0 && \ju \il = 0 \\
\label{e:ij1} & \iu \il \simeq \Id_{\cA} && \ju \jl \simeq \Id_{\cB}
\end{align}
and in addition,
\begin{align}
\label{e:ijC0} \Hom_{\cC}(\il A, \jl B) = 0 && \text{for all $A \in \ob
\cA$, $B \in \ob \cB$}
\end{align}

Using the adjunctions~\eqref{e:iadj} and~\eqref{e:jadj}, the
assumptions of~\eqref{e:ij0},\eqref{e:ij1} can be rephrased as follows.
First, the two conditions of~\eqref{e:ij0} are equivalent to each other
and each is equivalent to the condition
\begin{align}
\label{e:jiC0} \Hom_{\cC}(\jl B, \il A) = 0 && \text{for all $A \in \ob
\cA$, $B \in \ob \cB$}
\end{align}
Similarly, the conditions in~\eqref{e:ij1} are equivalent to the
requirement that $\il$ and $\jl$ are fully faithful functors, so that
one can think of $\cA$ and $\cB$ as embedded in $\cC$. Moreover,
from~\eqref{e:jiSES} and~\eqref{e:jiC0} we see that $(\cB, \cA)$ is a
\emph{torsion pair}~\cite[(I.2)]{HRS96} in $\cC$.

Observe also that~\eqref{e:ijC0} could be replaced with the assumption
that the functor $(\iu,\ju) : \cC \to \cA \times \cB$ is faithful.
Indeed, one implication follows from~\eqref{e:ij0} and~\eqref{e:ij1},
and the other follows using~\eqref{e:jiSES}.

\medskip

From now on assume that $(\iu, \il, \ju, \jl)$ form a gluing data.

\begin{lem} \label{l:iuproj}
If $P$ is a projective object of $\cC$, then $\iu P$ is projective in
$\cA$. Similarly, if $I$ an injective object of $\cC$, then $\ju I$ is
injective in $\cB$.
\end{lem}
\begin{proof}
A functor which is a left adjoint to an exact functor preserves
projectives, while a right adjoint to an exact functor preserves
injectives~\cite[Corollary~1.6]{FGR75}.
\end{proof}

The exact functors $\iu, \il, \ju, \jl$ give rise to triangulated
functors between the corresponding bounded derived categories. We use
the same notation for these derived functors:
\[
\xymatrix{\dA \ar@<1ex>[r]^{\il} & \dC \ar@<1ex>[l]^{\iu}
\ar@<1ex>[r]^{\ju} & \dB \ar@<1ex>[l]^{\jl}}
\]
Note that adjunctions and orthogonality relations analogous
to~\eqref{e:iadj}, \eqref{e:jadj}, \eqref{e:ij0}, \eqref{e:ij1} (but
not~\eqref{e:ijC0}) hold also for the derived functors. In particular,
$\il$ and $\jl$ are fully faithful.

\begin{defn}
An object $T$ in an abelian category $\cA$ is called \emph{rigid} if
$\Ext^n_{\cA}(T,T) = 0$ for all $n > 0$.
\end{defn}

\begin{prop}
\label{p:PIT} Let $P$ be a projective object of $\cC$ and $T_{\cB}$ be
a rigid object of $\cB$ such that $\Ext^n_{\cB}(\ju P, T_{\cB}) = 0$
for all $n > 0$. Consider the complex
\[
T = \ii P \oplus \jl T_{\cB} [1]
\]
in $\dC$. Then $\Hom_{\dC}(T, T[n]) = 0$ for $n \neq 0$ and
\[
\End_{\dC}(T) \simeq
\begin{pmatrix}
\End_{\cB}(T_{\cB}) & \Ext^1_{\cC}(\ii P, \jl T_{\cB}) \\
0 & \End_{\cA}(\iu P)
\end{pmatrix}
\]
is a triangular matrix ring.
\end{prop}
\begin{proof}
Since $T$ has two summands, the space $\Hom_{\dC}(T,T[n])$ decomposes
into four parts, which we now consider.

Since $\il$ is fully faithful and $\iu P$ is projective,
\begin{equation}
\label{e:iiP} \Hom_{\dC}(\ii P, \ii P[n]) \simeq \Hom_{\dA}(\iu P, \iu
P[n])
\end{equation}
vanishes for $n \neq 0$. Similarly, since $\jl$ is fully faithful and
$T_{\cB}$ is rigid,
\[
\Hom_{\dC}(\jl T_{\cB}, \jl T_{\cB}[n]) \simeq \Hom_{\dB}(T_{\cB},
T_{\cB}[n])
\]
vanishes for $n \neq 0$. Moreover, by orthogonality,
\[
\Hom_{\dC}(\jl T_{\cB}, \ii P[n]) = 0
\]
for all $n \in \bZ$.

It remains to consider $\Hom_{\dC}(\ii P, \jl T_{\cB}[n])$ and to prove
that it vanishes for $n \neq 1$. Using~\eqref{e:jiSES}, we obtain a
short exact sequence $0 \to \jj P \to P \to \ii P \to 0$ that induces a
long exact sequence, a fragment of which is shown below:
\begin{equation}
\label{e:PIles}
\begin{split}
\Hom_{\dB}(\ju P, T_{\cB}[n-1]) \simeq \Hom_{\dC}(\jj P[1], \jl T_{\cB}[n]) \to \\
\Hom_{\dC}(\ii P, \jl T_{\cB}[n]) \to \Hom_{\dC}(P, \jl T_{\cB}[n]) .
\end{split}
\end{equation}

Now observe that the right term vanishes for $n \neq 0$ since $P$ is
projective, and the left term of~\eqref{e:PIles} vanishes for $n \neq
1$ by our assumption on the vanishing of $\Ext^{\bullet}_{\cB}(\ju P,
T_{\cB})$. Therefore
\[
\Hom_{\dC}(\ii P,\jl T_{\cB}[n]) = 0
\]
for $n \neq 0, 1$. This holds also for $n=0$ by the
assumption~\eqref{e:ijC0}.

To complete the proof, note that $\Ext^1_{\cC}(\ii P, \jl T_{\cB})$ has
a natural structure of an $\End_{\cA}(\iu
P)$-$\End_{\cB}(T_{\cB})$-bimodule via the identifications
\begin{align*}
\End_{\cA}(\iu P) \simeq \End_{\cC}(\ii P) && \End_{\cB}(T_{\cB})
\simeq \End_{\cC}(\jl T_{\cB})
\end{align*}
\end{proof}

\begin{rem}
The assumptions in the proposition are always satisfied when $P$ is a
projective object of $\cC$ and $T_{\cB}$ is any injective object of
$\cB$.
\end{rem}

\begin{rem}
One can formulate an analogous statement for a rigid object $T_{\cA}$
of $\cA$ and an injective object $I$ of $\cC$.
\end{rem}

\section{Gluing in comma categories}
\label{sec:comma}

Let $\cA, \cB$ be categories and $F:\cA \to \cB$ a functor. The
\emph{comma category} with respect to the pair of functors $\cA
\xrightarrow{F} \cB \xleftarrow{\Id} \cB$ \cite[II.6]{MacLane98},
denoted by $(F \downarrow \Id)$, is the category $\cC$ whose objects
are triples $(A, B, f)$ where $A \in \ob \cA$, $B \in \ob \cB$ and $f :
FA \to B$ is a morphism (in $\cB$). The morphisms between objects
$(A,B,f)$ and $(A',B',f')$ are all pairs of morphisms $\alpha : A \to
A'$, $\beta : B \to B'$ such that the square
\begin{equation} \label{e:fAB}
\xymatrix{FA \ar[r]^{f} \ar[d]^{F\alpha} & B \ar[d]^{\beta} \\
FA' \ar[r]^{f'} & B'}
\end{equation}
commutes.

Assume in addition that $\cA$, $\cB$ are abelian and that $F: \cA \to
\cB$ is an additive, right exact functor. In this case, the comma
category $\cC$ is abelian \cite{FGR75}. Consider the functors
\[
\xymatrix{\cA \ar@<1ex>[r]^{\il} & \cC \ar@<1ex>[l]^{\iu}
\ar@<1ex>[r]^{\ju} & \cB \ar@<1ex>[l]^{\jl}}
\]
defined by
\begin{align*}
\il(A) = (A,0,0) && \il(\alpha) = (\alpha,0) & \qquad & \iu(A,B,f) = A
&& \iu(\alpha,\beta) = \alpha \\
\jl(B) = (0,B,0) && \jl(\beta) = (0,\beta) & \qquad & \ju(A,B,f) = B &&
\ju(\alpha,\beta) = \beta
\end{align*}
for objects $A \in \ob \cA$, $B \in \ob \cB$ and morphisms $\alpha$,
$\beta$.

\begin{lem}
The quadruple $(\iu,\il,\ju,\jl)$ is a gluing data.
\end{lem}
\begin{proof}
We need to verify the four properties of gluing data. The adjunction
follows by the commutativity of the diagrams
\begin{align*}
\xymatrix{FA \ar[r]^{f} \ar[d]^{F\alpha} & B \ar[d] \\
FA' \ar[r] & 0} &&
\xymatrix{0 \ar[r] \ar[d] & B \ar[d]^{\beta} \\
FA' \ar[r]^{f'} & B'}
\end{align*}
for $\alpha : A \to A'$ and $\beta : B \to B'$.

For exactness, note that kernels and images in $\cC$ can be computed
componentwise, that is, if $(\alpha,\beta) : (A,B,f) \to (A',B',f')$ is
a morphism in $\cC$, then $\ker (\alpha,\beta) = (\ker \alpha, \ker
\beta, f \vert_{F(\ker \alpha)})$ and similarly for the image. The
extension condition follows from
\[
0 \to (0,B,0) \xrightarrow{(0,1_B)} (A,B,f) \xrightarrow{(1_A,0)}
(A,0,0) \to 0
\]
and orthogonality is straightforward.
\end{proof}

One can use the special structure of the comma category $\cC$ to define
another pair of functors. Let $\ils : \cA \to \cC$ and $\jus : \cC \to
\cB$ be the functors defined by
\begin{align*}
&\ils(A) = (A,FA,1_{FA}) && \ils(\alpha) = (\alpha, F\alpha) \\
&\jus(A,B,f) = \coker f && \jus(\alpha,\beta) = \bar{\beta}
\end{align*}
where $\bar{\beta} : \coker f \to \coker f'$ is induced from $\beta$.

\begin{lem} \label{l:ils}
\label{l:ijLeft} $\ils$ is a left adjoint of $\iu$, $\jus$ is a left
adjoint of $\jl$, and
\begin{align*}
\iu \ils \simeq \Id_{\cA} && \ju \ils = F && \jus \ils = \jus \il = 0
&& \jus \jl \simeq \Id_{\cB}
\end{align*}
\end{lem}
\begin{proof}
The adjunctions follow by considering the commutative diagrams
\begin{align*}
\xymatrix{FA \ar[r]^{1_{FA}} \ar[d]^{F\alpha} & FA \ar[d]^{\beta = f' \circ F\alpha} \\
FA' \ar[r]^{f'} & B'} &&
\xymatrix{FA \ar[r]^{f} \ar[d] & B \ar[d]^{\beta} \\
0 \ar[r] & B'}
\end{align*}
and noting that the commutativity of the right diagram implies that
$\beta$ factors uniquely through $\coker f$. The other relations are
straightforward.
\end{proof}

\begin{rem}
The diagram
\[
\xymatrix{\cA \times \cB \ar@<1ex>[r]^{(\ils,\jl)} & \cC
\ar@<1ex>[l]^{(\iu,\ju)} \ar@<1ex>[r]^{(\iu,\jus)} & \cA \times \cB
\ar@<1ex>[l]^{(\il,\jl)}}
\]
is a special case of the one in~\cite[p.~7]{FGR75}, viewing $\cC$ as a
trivial extension of $\cA \times \cB$.
\end{rem}

\begin{prop} \label{p:PITcomma}
Let $P$ be a projective object of $\cA$ and $T_{\cB}$ be a rigid object
of $\cB$ such that $\Ext^n_{\cB}(FP, T_{\cB}) = 0$ for all $n > 0$.
Assume that $FP \in \ob(\cB)$ has a projective resolution in $\cB$ and
consider $T = (P, 0, 0) \oplus (0, T_{\cB}, 0)[1] \in \dC$. Then
$\Hom_{\dC}(T, T[n]) = 0$ for $n \neq 0$ and
\[
\End_{\dC}(T) \simeq
\begin{pmatrix}
\End_{\cB}(T_{\cB}) & \Hom_{\cB}(FP, T_{\cB}) \\
0 & \End_{\cA}(P)
\end{pmatrix},
\]
where the bimodule structure on $\Hom_{\cB}(FP, T_{\cB})$ is given by
left composition with $\End_{\cB}(T_{\cB})$ and right composition with
$\End_{\cB}(FP)$ through the natural map $\End_{\cA}(P) \to
\End_{\cB}(FP)$.
\end{prop}
\begin{proof}
Since $\ils$ is a left adjoint of an exact functor, it takes projective
objects of $\cA$ to projective objects of $\cC$. Hence $\ils P = (P,
FP, 1_{FP})$ is projective and we can apply Proposition~\ref{p:PIT} for
$\ils P$ and $T_{\cB}$. As $(P,0,0) = \ii \ils P$ and
$(0,T_{\cB},0)=\jl T_{\cB}$, we only need to show the isomorphism
\[
\Ext^1_{\cC}\left((P,0,0),(0,T_{\cB},0)\right) \simeq \Hom_{\cB}(FP,
T_{\cB})
\]
as $\End_{\cA}(P)$-$\End_{\cB}(T_{\cB})$-bimodules.

Indeed, let
\begin{equation} \label{e:projresFP}
\dots \to Q^2 \to Q^1 \to FP \to 0
\end{equation}
be a projective resolution of $FP$. Then $(P,0,0)$ is quasi-isomorphic
to the complex
\[
\dots \to \jl Q^2 \to \jl Q^1 \to \ils P \to 0 \to \dots
\]
whose terms are projective since $\jl$ is a left adjoint of an exact
functor. Therefore $\Ext^1_{\cC}((P,0,0),(0,T_{\cB},0))$ can be
identified with the morphisms, up to homotopy, between the complexes
\begin{equation} \label{e:PprojI}
\xymatrix{
\dots \ar[r] & \jl Q^2 \ar[r] & \jl Q^1 \ar[r] &
\ils P \ar[r] & 0 \ar[r] & \dots \\
\dots \ar[r] & 0 \ar[r] & \jl T_{\cB} \ar[r] &
0 \ar[r] & 0 \ar[r] & \dots
}
\end{equation}

By Lemma~\ref{l:ils}, $\Hom_{\cC}(\ils P, \jl T_{\cB}) = \Hom_{\cA}(P,
\iu \jl T_{\cB}) = 0$, thus any homotopy between these complexes
vanishes, and the morphism space equals $\ker(\Hom_{\cC}(\jl Q^1, \jl
T_{\cB}) \to \Hom_{\cC}(\jl Q^2, \jl T_{\cB}))$. Using the fact that
$\jl$ is fully faithful and applying the functor $\Hom_{\cB}(-,
T_{\cB})$ on the exact sequence~\eqref{e:projresFP}, we get that the
morphism space equals $\Hom_{\cB}(FP, T_{\cB})$, as desired.

Under this identification, the left action of $\End_{\cB}(T_{\cB})
\simeq \End_{\cC}(\jl T_{\cB})$ is given by left composition. As for
the right action of $\End_{\cA}(P)$, observe that any $\alpha \in
\End_{\cA}(P)$ extends uniquely to an endomorphism in the homotopy
category
\[
\xymatrix{(0, FP, 0) \ar[r]^{(0,1)} \ar[d]^{(0, F\alpha)} &
(P,FP, 1_{FP}) \ar[d]^{(\alpha, F\alpha)} \\
(0, FP, 0) \ar[r]^{(0,1)} & (P, FP, 1_{FP})}
\]
which determines a unique endomorphism, in the homotopy category, of
the top complex of~\eqref{e:PprojI}.
\end{proof}

\begin{rem} \label{r:rightadj}
When the functor $F : \cA \to \cB$ admits a right adjoint $G : \cB \to
\cA$, the comma category $(F \downarrow \Id)$ is equivalent to the
comma category $(\Id \downarrow G)$ corresponding to the pair $\cA
\xrightarrow{\Id} \cA \xleftarrow{G} \cB$. In this case, one can define
also a right adjoint $\ius$ of $\il$ and a right adjoint $\jls$ of
$\ju$, and we end up with the eight functors $(\ils, \iu, \il, \ius)$
and $(\jus, \jl, \ju, \jls)$. The bimodule $\Hom_{\cB}(FP, T_{\cB})$ in
Proposition~\ref{p:PITcomma} can then be identified with $\Hom_{\cA}(P,
G T_{\cB})$.
\end{rem}

\section{Application to triangular matrix rings}
\label{sec:alg}

\subsection{Triangular matrix rings}
Let $R$ and $S$ be rings, and let $_RM_S$ be an $S$-$R$-bimodule. Let
$\gL$ be the \emph{triangular matrix ring}
\begin{equation} \label{e:Ltriang}
\gL = \begin{pmatrix} R & M \\ 0 & S \end{pmatrix} = \left\{
\begin{pmatrix} r & m \\ 0 & s \end{pmatrix} \,:\, r \in R, s \in S, m
\in M \right\}
\end{equation}
where the ring structure is induced by the ordinary matrix operations.

For a ring $R$, denote the category of right $R$-modules by $\Mod R$.
The functor $- \otimes M : \Mod R \to \Mod S$ is additive and right
exact, hence the corresponding comma category $(- \otimes M \downarrow
\Id_{\Mod R})$ is abelian.

\begin{lem}[\protect{\cite[III.2]{ARS95}}] \label{l:RSL}
The category $\Mod \gL$ is equivalent to the comma category $(- \otimes
M \downarrow \Id_{\Mod R})$.
\end{lem}
\begin{proof}
One verifies that by sending a triple $(X_R, Y_S, f : X \otimes M \to
Y)$ to the $\gL$-module $X \oplus Y$ defined by
\begin{equation} \label{e:RSL}
\begin{pmatrix}
  x & y
\end{pmatrix}
\begin{pmatrix}
  r & m \\ 0 & s
\end{pmatrix}
=
\begin{pmatrix}
  xr & f(x \otimes m) + ys
\end{pmatrix}
\end{equation}
and sending a morphism $(\alpha,\beta) : (X,Y,f) \to (X',Y',f')$ to
$\alpha \oplus \beta : X \oplus Y \to X' \oplus Y'$, we get a functor
$(- \otimes M \downarrow \Id_{\Mod R}) \to \Mod \gL$ which is an
equivalence of categories.
\end{proof}

\begin{cor}
\label{c:glueMod} There exists gluing data
\[
\xymatrix{\Mod R \ar@<1ex>[r]^{\il} & \Mod \gL \ar@<1ex>[l]^{\iu}
\ar@<1ex>[r]^{\ju} & \Mod S. \ar@<1ex>[l]^{\jl}}
\]
\end{cor}

The functors occurring in Corollary~\ref{c:glueMod} can be described
explicitly. Let
\begin{align*}
e_R = \begin{pmatrix} 1 & 0 \\ 0 & 0 \end{pmatrix} &,& e_S =
\begin{pmatrix} 0 & 0 \\ 0 & 1 \end{pmatrix}.
\end{align*}
Using~\eqref{e:RSL}, observe that for a $\gL$-module $Z_{\gL}$,
\begin{align} \label{e:iuju}
(\iu Z)_R = Z e_R && (\ju Z)_S = Z e_S
\end{align}
where $r$ acts on $\iu Z$ via $\left(\begin{smallmatrix} r & 0 \\ 0 & 0
\end{smallmatrix}\right)$ and $s$ acts on $\ju Z$ via $\left(\begin{smallmatrix}
0 & 0 \\ 0 & s \end{smallmatrix}\right)$. The morphism $(\iu Z) \otimes
M \to \ju Z$ is obtained by considering the actions of
$\left(\begin{smallmatrix} 0 & m \\ 0 & 0\end{smallmatrix}\right)$, $m
\in M$, and the map $Z \mapsto (\iu Z, \ju Z, (\iu Z) \otimes M \to \ju
Z)$ defines a functor which is an inverse to the equivalence of
categories constructed in the proof of Lemma~\ref{l:RSL}.

Conversely, for an $R$-module $X_R$ and $S$-module $Y_S$, we have $(\il
X)_{\gL} = X$ and $(\jl Y)_{\gL} = Y$ where $\left(\begin{smallmatrix}
r & m \\ 0 & s
\end{smallmatrix}\right)$ acts on $X$ via $r$ and on $Y$ via $s$.

\begin{lem}
The image of $\gL_{\gL}$ in the comma category equals $(R, M, 1_M)
\oplus (0, S, 0)$.
\end{lem}
\begin{proof}
Use~\eqref{e:iuju} and
\begin{align*}
\gL \begin{pmatrix} 1 & 0 \\ 0 & 0 \end{pmatrix} = \begin{pmatrix} R &
0 \\ 0 & 0 \end{pmatrix} &,& \gL \begin{pmatrix} 0 & 0 \\ 0 & 1
\end{pmatrix} = \begin{pmatrix} 0 & M \\ 0 & S \end{pmatrix}.
\end{align*}
\end{proof}

\begin{rem}
Since $- \otimes M$ admits a right adjoint $\Hom(M, -)$, we are in the
situation of Remark~\ref{r:rightadj} and there are eight functors
$(\ils, \iu, \il, \ius)$ and $(\jus, \jl, \ju, \jls)$. For the
convenience of the reader, we now describe them as standard functors
$\otimes$ and $\Hom$ involving idempotents, see
also~\cite[Section~2]{CPS88b} and~\cite[Proposition~2.17]{Miyachi03}.

If $A$ is a ring and $e \in A$ is an idempotent, the functor
\[
\Hom_{A}(eA, -) = - \otimes_{A} Ae : \Mod A \to \Mod eAe
\]
admits a left adjoint $- \otimes_{eAe} eA$ and a right adjoint
$\Hom_{eAe}(Ae, -)$. By taking $A = \gL$ and $e = e_R$ we get the three
functors $(\ils, \iu, \il)$. Similarly, $e=e_S$ gives $(\jl, \ju,
\jls)$.

In addition, the natural inclusion functor
\[
\Hom_{A/AeA}(A/AeA, -) = - \otimes_{A/AeA} A/AeA : \Mod A/AeA \to \Mod
A
\]
admits a left adjoint $- \otimes_A A/AeA$ and a right adjoint
$\Hom_A(A/AeA, -)$. By taking $A = \gL$ and $e = e_R$, observing that
$e_S \gL e_R = 0$, we get the three functors $(\jus, \jl, \ju)$. For $e
= e_S$ we get $(\iu, \il, \ius)$.
\end{rem}

\subsection{The main theorem}
For a ring $\gL$, denote by $\dL$ the bounded derived category of $\Mod
\gL$, and by $\per \gL$ its full subcategory of complexes
quasi-isomorphic to \emph{perfect complexes}, that is, bounded
complexes of finitely generated projective $\gL$-modules.

For a complex $T \in \dL$, denote by $\langle \add T \rangle$ the
smallest full triangulated subcategory of $\dL$ containing $T$ and
closed under forming direct summands. Recall that $T$ is a
\emph{tilting complex} if $\langle \add T \rangle = \per \gL$ and
$\Hom_{\dL}(T, T[n]) = 0$ for all integers $n \neq 0$. If, moreover,
$H^n(T) = 0$ for all $n \neq 0$, we call $T$ a \emph{tilting module}
and identify it with the module $H^0(T)$.

\begin{thm} \label{t:RSM}
Let $R, S$ be rings and $T_S$ a tilting $S$-module. Let $_RM_S$ be an
$S$-$R$-bimodule such that as an $S$-module, $M_S \in \per S$ and
$\Ext^n_S(M_S, T_S) = 0$ for all $n > 0$. Then the triangular matrix
rings
\[
\gL =
\begin{pmatrix}
R & M \\
0 & S \\
\end{pmatrix}
\quad \text{and} \quad \widetilde{\gL} =
\begin{pmatrix}
\End_S(T_S) & \Hom_S(M, T_S) \\
0 & R \\
\end{pmatrix}
\]
are derived equivalent.
\end{thm}
\begin{proof}
For simplicity, we shall identify $\Mod \gL$ with the corresponding
comma category. We will show that $T = (R, 0, 0) \oplus (0, T_S, 0)[1]$
is a tilting complex in $\dL$ whose endomorphism ring equals
$\widetilde{\gL}$.

Applying Proposition~\ref{p:PITcomma} for the projective $R$-module $R$
and the rigid $S$-module $T_S$, noting that $FR = M_S$ and
$\Ext^n_S(M_S, T_S) = 0$ for $n > 0$, we see that $\Hom_{\dL}(T, T[n])
= 0$ for all $n \neq 0$ and moreover $\End_{\dL}(T) \simeq
\widetilde{\gL}$.

It remains to show that $\langle \add T \rangle = \per \gL$. First, we
show that $T \in \per \gL$. Observe that $\jl(\per S) \subseteq \per
\gL$, since $\jl$ is an exact functor which takes projectives to
projectives and $\jl S = (0, S, 0)$ is a direct summand of $\gL$. Hence
in the short exact sequence
\begin{equation} \label{e:RMperf}
0 \to (0, M, 0) \to (R, M, 1_M) \to (R, 0, 0) \to 0 ,
\end{equation}
we have that $(0, M, 0) \in \per \gL$ by the assumption that $M_S \in
\per S$, and $(R, M, 1_M) \in \per \gL$ as a direct summand of $\gL$.
Therefore $(R, 0, 0) \in \per \gL$. In addition, $(0, T_S, 0) \in \per
\gL$ by the assumption $T_S \in \per S$, hence $T$ is isomorphic in
$\dL$ to a perfect complex.

Second, in order to prove that $\langle \add T \rangle = \per \gL$ it
is enough to show that $\gL \in \langle \add T \rangle$. Indeed, since
$(0, T_S, 0)[1]$ is a summand of $T$, by the exactness of $\jl$ and our
assumption that $\langle \add T_S \rangle = \per S$, we have that $(0,
S, 0) \in \langle \add T \rangle$ and $(0, M, 0) \in \langle \add T
\rangle$. Since $(R, 0, 0)$ is a summand of $T$, by invoking again the
short exact sequence~\eqref{e:RMperf} we see that $(R, M, 1_M) \in
\langle \add T \rangle$, hence $\gL \in \langle \add T \rangle$.

Therefore $T$ is a tilting complex in $\dL$, and by \cite{Rickard89}
(see also~\cite[(1.4)]{Keller98}), the rings $\gL$ and $\widetilde{\gL}
\simeq \End_{\dL}(T)$ are derived equivalent.
\end{proof}

\begin{rem}
The assumption that $T_S$ is a tilting module implies that the rings
$S$ and $\End_S(T_S)$ are derived equivalent.
\end{rem}

\begin{rem} \label{r:RSMinj}
When the tilting module $T_S$ is also injective, it is enough to assume
that $M_S \in \per S$.
\end{rem}

\subsection{Applications}

Let $S$ be an Artin algebra over an Artinian commutative ring $k$, and
let $\modf S$ be the category of finitely generated right $S$-modules.
Let $D: \Mod S \to \Mod S^{op}$ be the functor defined by $D =
\Hom_k(-, J)$, where $J$ is an injective envelope of the direct sum of
all the non-isomorphic simple modules of $k$. Recall that $D$ restricts
to a duality $D : \modf S \to \modf S^{op}$. Applying it on the
bimodule $_SS_S$, we get the bimodule $_SDS_S = \Hom_k(S, J)$.

\begin{lem} \label{l:DM}
Let $R$ be a ring and $_RM_S$ a bimodule. Then $_SDM_R \simeq
\Hom_S(_RM_S, {_SDS_S})$ as $R$-$S$-bimodules.
\end{lem}
\begin{proof}
By the standard adjunctions,
\[
\Hom_S(_RM_S, \Hom_k(_SS_S, J)) \simeq \Hom_k(_RM_S \otimes _SS_S, J) =
\Hom_k(_RM_S, J).
\]
\end{proof}

It follows that $D = \Hom_S(-, DS_S)$, hence $DS_S$ is an injective
object in $\Mod S$. We denote by $\gldim S$ the global dimension of
$\modf S$.

\begin{thm} \label{t:RSMartin}
Let $R$ be a ring, $S$ an Artin algebra with $\gldim S < \infty$ and
$_RM_S$ an $S$-$R$-bimodule which is finitely generated as an
$S$-module. Then the triangular matrix rings
\[
\gL =
\begin{pmatrix}
R & M \\
0 & S \\
\end{pmatrix}
\quad \text{and} \quad \widetilde{\gL} =
\begin{pmatrix}
S & DM \\
0 & R \\
\end{pmatrix}
\]
are derived equivalent, where $D$ is the duality on $\modf S$.
\end{thm}
\begin{proof}
The module $DS_S$ is injective in $\Mod S$ and any module in $\modf S$
has an injective resolution with terms that are summands of finite
direct sums of $DS$. Since $\gldim S < \infty$, such a resolution is
finite, hence $\langle \add DS \rangle = \per S$ and $M \in \per S$ for
any $M \in \modf S$.

Therefore the assumptions of Theorem~\ref{t:RSM} are satisfied for $T_S
= DS$ (see also Remark~\ref{r:RSMinj}), and it remains to show that
$\End_S(T_S) = S$ and $\Hom_S(M, T_S) \simeq {_SDM_R}$ (as bimodules).
This follows by the Lemma~\ref{l:DM} applied for the bimodules $_SDS_S$
and $_RM_S$.
\end{proof}

\begin{rem} \label{r:RSMartinalg}
Under the assumptions of Theorem~\ref{t:RSMartin}, when $R$ is also an
Artin $k$-algebra and $k$ acts centrally on $M$, the rings $\gL$ and
$\widetilde{\gL}$ are Artin algebras and the derived equivalence in the
theorem implies that $\cD^b(\modf \gL) \simeq \cD^b(\modf
\widetilde{\gL})$.

Moreover, by using the duality $D$, one sees that
Theorem~\ref{t:RSMartin} is true for two Artin algebras $R$ and $S$ and
a bimodule $_RM_S$ on which $k$ acts centrally under the weaker
assumptions that $M$ is finitely generated over $k$ and at least one of
$\gldim R$, $\gldim S$ is finite.
\end{rem}

\medskip

By taking $T_S = S$ in Theorem~\ref{t:RSM}, we get the following
corollary.

\begin{cor} \label{c:RSMproj}
Let $R$, $S$ be rings and $_RM_S$ an $S$-$R$-bimodule such that as an
$S$-module, $M_S \in \per S$ and $\Ext^n_S(M_S, S) = 0$ for all $n >
0$. Then the triangular matrix rings
\[
\gL =
\begin{pmatrix}
R & M \\
0 & S \\
\end{pmatrix}
\quad \text{and} \quad \widetilde{\gL} =
\begin{pmatrix}
S & \Hom_S(M, S) \\
0 & R \\
\end{pmatrix}
\]
are derived equivalent.
\end{cor}

\begin{rem}
The conditions of Corollary~\ref{c:RSMproj} hold when the ring $S$ is
\emph{self-injective}, that is, $S$ is injective as a (right) module
over itself, and $_RM_S$ is finitely generated projective as an
$S$-module. In particular, this applies when $S$ is a semi-simple ring
and $M$ is finitely generated as an $S$-module.
\end{rem}

\begin{rem}
Recall that for a ring $R$, a division ring $S$ and a bimodule $_SN_R$
which is finite dimensional as a left $S$-vector space, the
\emph{one-point extension} $R[N]$ and the \emph{one-point coextension}
$[N]R$ of $R$ by $N$ are defined as the triangular matrix rings
\begin{align*}
R[N] =
\begin{pmatrix}
S & _SN_R \\ 0 & R
\end{pmatrix}
&& [N]R =
\begin{pmatrix}
R & _RDN_S \\ 0 & S
\end{pmatrix}.
\end{align*}
where $D = \Hom_S(-, S)$ is the duality on $\modf S$. By taking $M =
DN$ in the preceding remark, we see that the rings $R[N]$ and $[N]R$
are derived equivalent. Compare this with the construction of
``reflection with respect to an idempotent''
in~\cite{TachikawaWakamatsu86}.
\end{rem}

\section{Concluding remarks}
\label{sec:remarks}

\subsection{Repetitive algebras}
In the specific case of Artin algebras, another approach to the
connection between a triangular matrix algebra $\gL$ and its mate
$\widetilde{\gL}$ involves the use of repetitive algebras, as outlined
below.

Let $\gL$ be an Artin algebra over a commutative Artinian ring $k$ and
let $D: \modf k \to \modf k$ be the duality. Recall that the
\emph{repetitive algebra} $\widehat{\gL}$ of $\gL$, introduced
in~\cite{HughesWaschbusch83}, is the algebra (without unit) of matrices
of the form
\[
\widehat{\gL} =
\begin{pmatrix}
\ddots & D\gL_{i-1} & 0 \\
0 & \gL_i & D\gL_i & 0 \\
& 0 & \gL_{i+1} & D\gL_{i+1} \\
& & 0 & \ddots
\end{pmatrix}
\]
where $\gL_i = \gL$, $D\gL_i = D \gL$ for $i \in \bZ$, and only finite
number of entries are nonzero. The multiplication is defined by the
canonical maps $\gL \otimes_{\gL} D\gL \to D\gL$, $D\gL \otimes_{\gL}
\gL \to D\gL$ induced by the bimodule structure on $D\gL$, and the zero
map $D\gL \otimes_{\gL} D\gL \to 0$.

When $\gL$ is a triangular matrix algebra, one can write
\begin{align*}
\gL = \begin{pmatrix} R & M \\ 0 & S \end{pmatrix} && D\gL
=\begin{pmatrix} DR & 0 \\ DM & DS \end{pmatrix}
\end{align*}
and a direct calculation shows that the maps $\gL \otimes D\gL \to
D\gL$ and $D\gL \otimes \gL \to D\gL$ are given by multiplication of
the above matrices, under the convention that $M \otimes_S DS \to 0$
and $DR \otimes_R M \to 0$.

As for the mate $\widetilde{\gL}$, we have
\begin{align*}
\widetilde{\gL} = \begin{pmatrix} S & DM \\ 0 & R \end{pmatrix} && D
\widetilde{\gL} =\begin{pmatrix} DS & 0 \\ M & DR \end{pmatrix},
\end{align*}
therefore the repetitive algebras of $\gL$ and its mate
$\widetilde{\gL}$ have the form
\begin{align*}
\widehat{\gL} = \left(
\begin{smallmatrix}
\ddots & DM & DS \\
& R & M & DR \\
& & S & DM & DS \\
& & & R & M & DR \\
& & & & S & DM & DS\\
& & & & & R & M \\
& & & & & & \ddots
\end{smallmatrix}
\right) &,& \widehat{\widetilde{\gL}} = \left(
\begin{smallmatrix}
\ddots & M & DR \\
& S & DM & DS \\
& & R & M & DR \\
& & & S & DM & DS \\
& & & & R & M & DR\\
& & & & & S & DM \\
& & & & & & \ddots
\end{smallmatrix}
\right)
\end{align*}
and are thus clearly seen to be isomorphic.

When $k$ is a field and \emph{both} algebras $R$ and $S$ have finite
global dimension, this can be combined with Happel's
Theorem~\cite[(II.4.9)]{Happel88} to deduce that $\gL$ and its mate
$\widetilde{\gL}$ are derived equivalent.

Note, however, that for the derived equivalence between $\gL$ and
$\widetilde{\gL}$ to hold, it is enough to assume that only \emph{one}
of $R$, $S$ has finite global dimension (see
Remark~\ref{r:RSMartinalg}).

Moreover, while the repetitive algebras of $\gL$ and $\widetilde{\gL}$
are always isomorphic, in the case where none of $R$, $S$ have finite
global dimension, the algebras $\gL$ and $\widetilde{\gL}$ may not be
derived equivalent, see Example~\ref{ex:matenoneq} below.

\subsection{Grothendieck groups}

In this subsection, $k$ denotes an algebraically closed field. Let
$\gL$ be a finite dimensional $k$-algebra and let $P_1,\dots,P_n$ be a
complete collection of the non-isomorphic indecomposable projectives in
$\modf \gL$. The \emph{Cartan matrix} of $\gL$ is the $n \times n$
integer matrix defined by $C_{ij} = \dim_k \Hom(P_i, P_j)$.

The Grothendieck group $K_0(\per \gL)$ of the triangulated category
$\per \gL$ can be viewed as a free abelian group on the generators
$[P_1],\dots,[P_n]$, and the Euler form
\[
\langle K, L \rangle = \sum_{r \in \bZ} (-1)^r \dim_k \Hom_{\dL}(K,
L[r])
\]
on $\per \gL$ induces a bilinear form on $K_0(\per \gL)$ whose matrix
with respect to that basis equals the Cartan matrix.

It is well known that a derived equivalence of two algebras $\gL$ and
$\gL'$ induces an equivalence of the triangulated categories $\per \gL$
and $\per \gL'$, and hence an isometry of their Grothendieck groups
preserving the Euler forms. We now consider the consequences of the
derived equivalence of Theorem~\ref{t:RSMartin} (when $R$ and $S$ are
finite dimensional $k$-algebras) for the corresponding Grothen\-dieck
groups.

For simplicity, assume that $\gL$ is basic. In this case, there exist
primitive orthogonal idempotents $\{e_1,\dots,e_n\}$ in $\gL$ such that
$P_i \simeq e_i \gL$ for $1 \leq i \leq n$. Therefore by the
isomorphisms $\Hom_{\gL}(e_i \gL, N) \simeq N e_i$ of $k$-spaces for
any $\gL$-module $N_{\gL}$, we get that $C_{ij} = \dim_k e_j \gL e_i$.

\begin{lem} \label{l:Cartan}
Let $R$, $S$ be basic, finite dimensional $k$-algebras, and let $_RM_S$
be a finite dimensional $S$-$R$-bimodule. Then the Cartan matrix
$C_{\gL}$ of the corresponding triangular matrix algebra $\gL$ is of
the form
\[
C_{\gL} = \begin{pmatrix} C_R & 0 \\ C_M & C_S
\end{pmatrix}
\]
where $C_R$, $C_S$ are the Cartan matrices of $R$, $S$.
\end{lem}
\begin{proof}
Let $e_1,\dots,e_n$ and $f_1,\dots,f_m$ be complete sets of primitive
orthogonal idempotents in $R$ and in $S$. Let $\bar{e}_i = e_i
\left(\begin{smallmatrix} 1 & 0 \\ 0 & 0 \end{smallmatrix}\right)$ and
$\bar{f}_j = f_j \left(\begin{smallmatrix} 0 & 0 \\ 0 & 1
\end{smallmatrix}\right)$. Then
$\bar{e}_1,\dots,\bar{e}_n,\bar{f}_1,\dots,\bar{f}_m$ is a complete set
of primitive orthogonal idempotents of $\gL$ and the result follows by
computing the dimensions of $\bar{e}_i \gL \bar{e}_{i'}$, $\bar{e}_i
\gL \bar{f}_j$, $\bar{f}_j \gL \bar{e}_i$ and $\bar{f}_j \gL
\bar{f}_{j'}$. In particular, $(C_M)_{ji} = \dim_k e_i M f_j$.
\end{proof}

Since $\dim_k f_j DM e_i = \dim_k e_i M f_j$, we get by
Lemma~\ref{l:Cartan} that the Cartan matrices of $\gL$ and its mate
$\widetilde{\gL}$ are
\begin{align*}
C_{\gL} =
\begin{pmatrix}
C_R & 0 \\ C_M & C_S
\end{pmatrix}
&& C_{\widetilde{\gL}} =
\begin{pmatrix}
C_S & 0 \\ C_M^t & C_R
\end{pmatrix}.
\end{align*}

When at least one of $R$ and $S$ has finite global dimension, the
derived equivalence of Theorem~\ref{t:RSMartin} implies that $C_{\gL}$
and $C_{\widetilde{\gL}}$ represent the same bilinear form, hence they
are congruent over $\bZ$, that is, there exists an invertible matrix
$P$ over $\bZ$ such that $P^t C_{\gL} P = C_{\widetilde{\gL}}$.

One can also show this congruence directly at the level of matrices, as follows.

\begin{lem} \label{l:congmat}
Let $K$ be a commutative ring. Let $A \in M_{n \times n}(K)$ be a
square matrix, $B \in \GL_m(K)$ an invertible square matrix and $C \in
M_{m \times n}(K)$. Then there exists $P \in \GL_{n+m}(K)$ such that
\[
P^t \begin{pmatrix} A & 0 \\ C & B \end{pmatrix} P = \begin{pmatrix} B&
0 \\ C^t & A \end{pmatrix}
\]
\end{lem}
\begin{proof}
Take $P = \begin{pmatrix} 0 & I_n \\ -B^{-1} B^t & -B^{-1} C
\end{pmatrix}$. Then
\begin{align*}
P^t \begin{pmatrix} A & 0 \\ C & B \end{pmatrix} P &=
\begin{pmatrix}
0 & -B B^{-t} \\ I_n & -C^t B^{-t}
\end{pmatrix}
\begin{pmatrix}
A & 0 \\ C & B
\end{pmatrix}
\begin{pmatrix}
0 & I_n \\ -B^{-1} B^t & -B^{-1} C
\end{pmatrix}
\\
&=
\begin{pmatrix}
0 & -B B^{-t} \\ I_n & -C^t B^{-t}
\end{pmatrix}
\begin{pmatrix}
0 & A \\ -B^t & 0
\end{pmatrix}
=
\begin{pmatrix}
B & 0 \\ C^t & A
\end{pmatrix}
\end{align*}
\end{proof}

Note that one could also take $P = \begin{pmatrix} -A^{-t} C^t &
-A^{-t} A \\ I_m & 0 \end{pmatrix}$, hence it is enough to assume that
at least one of $A$ and $B$ is invertible.

The conclusion of the lemma is false if one does not assume that at
least one of the matrices $A$, $B$ is invertible over $K$. This can be
used to construct triplets consisting of two finite dimensional
algebras $R$, $S$ (necessarily of infinite global dimension) and a
bimodule $M$ such that the triangular matrix algebra $\gL$ and its mate
$\widetilde{\gL}$ are not derived equivalent.

\begin{exmp} \label{ex:matenoneq}
Let $R = k[x]/(x^2)$, $S = k[y]/(y^3)$ and $M=k$ with $x$ and $y$
acting on $k$ as zero. Then the triangular matrix algebras
\begin{align*}
\gL =
\begin{pmatrix}
k[x]/(x^2) & k \\ 0 & k[y]/(y^3)
\end{pmatrix}
&& \widetilde{\gL} =
\begin{pmatrix}
k[y]/(y^3) & k \\ 0 & k[x]/(x^2)
\end{pmatrix}
\end{align*}
are not derived equivalent, since one can verify that their Cartan
matrices
\begin{align*}
C_{\gL} =
\begin{pmatrix}
2 & 0 \\ 1 & 3
\end{pmatrix}
&& C_{\widetilde{\gL}} =
\begin{pmatrix}
3 & 0 \\ 1 & 2
\end{pmatrix}
\end{align*}
are not congruent over $\bZ$. Note that despite the fact that $R$ and
$S$ are self-injective, Corollary~\ref{c:RSMproj} cannot be used since
$M$ does not have a finite projective resolution.
\end{exmp}

\subsection{Trivial extensions}
Triangular matrix rings are special cases of trivial
extensions~\cite[p.~78]{ARS95}. Indeed, if $R$, $S$ are rings and
$_RM_S$ is a bimodule, the corresponding triangular matrix ring is
isomorphic to the trivial extension $A \ltimes M$ where $A = R \times
S$ and $M$ is equipped with an $A$-bimodule structure via $(r,s)m=rm$
and $m(r,s)=ms$.

We remark that even when $A$ is a finite dimensional $k$-algebra of
finite global dimension and $M$ is a finite dimensional $A$-bimodule,
the trivial extension algebras $A \ltimes M$ and $A \ltimes DM$ are
generally not derived equivalent, so that the derived equivalence in
Theorem~\ref{t:RSMartin} is a special feature of triangular matrix
rings.

\begin{exmp}
Let $A = kQ$ where $Q$ is the quiver $\xymatrix{{\bullet_1} \ar[r]&
{\bullet_2}}$ and let $M$ be the $kQ$-bimodule corresponding to the
following commutative diagram of vector spaces
\[
\xymatrix@=1pc{& 0_{(1,1)} \ar[dr] \\
0_{(2,1)} \ar[ur] \ar[dr] & & {k}_{(1,2)} \\
& 0_{(2,2)} \ar[ur]}
\]
Then $A \ltimes M$ is the path algebra of the quiver ${\bullet}
\rightrightarrows {\bullet}$ while $A \ltimes DM$ is the path algebra
of ${\bullet} \rightleftarrows {\bullet}$ modulo the compositions of
the arrows being zero. These two algebras are not derived equivalent
since $\gldim (A \ltimes M) = 1$ while $\gldim (A \ltimes DM) =
\infty$.
\end{exmp}


\end{document}